\documentclass[10pt,twocolumn]{article}
\usepackage{makeidx}
\usepackage[french, english]{babel}
\usepackage[latin1]{inputenc}
\usepackage{amsfonts}
\usepackage{bm}
\usepackage{amssymb}
\usepackage{latexsym}
\usepackage{amsmath}
\usepackage{amscd}
\usepackage{amsthm}
\usepackage[colorlinks=true]{hyperref}
\usepackage{rotating}
\usepackage{graphicx}
\usepackage{pifont}
\usepackage{ulem}
\usepackage{curves}
\usepackage{fancyhdr}
\usepackage{epsfig}
\usepackage{pstricks}
\usepackage{subfigure}
\usepackage{epic}
\usepackage{alltt}

\newcommand{\E}{\mathbb{E}}
\renewcommand{\P}{\mathbb{P}}
\newcommand{\om}{\omega}

\newcommand{\ind}{\mathbf{1}}

\newcommand{\Co}{\mathcal{C}}

\newcommand{\N}{\mathbb{N}}
\newcommand{\Z}{\mathbb{Z}}
\newcommand{\R}{\mathbb{R}}

\newcommand{\etal}{\textit{et al.} }
\newcommand{\eg}{\textit{e.g.} }
\newcommand{\eps}{\varepsilon}
\newcommand{\e}{\mathrm{e}}

\renewcommand{\a}{\alpha}
\renewcommand{\b}{\beta}
\renewcommand{\l}{\lambda}
\renewcommand{\L}{\Lambda}
\newcommand{\vm}[1]{\E\lambda(#1)}
\newcommand{\et}{\quad\text{and}\quad}
\newcommand{\cube}{{[0,1]^d}}
\newcommand{\dimb}{\overline{\mathrm{dim}_\mathrm{box}}}
\renewcommand{\d}{\mathrm{d}}

\newcommand{\Om}{\Omega}

\newcommand{\dd}{\boldsymbol d}

\newcommand{\1}{{\bf 1}}

\theoremstyle{plain}
\newtheorem{theorem}{Theorem}[section]
\newtheorem{proposition}[theorem]{Proposition}
\newtheorem{lemma}[theorem]{Lemma}

\newtheorem{definition}[theorem]{Definition}

\newtheorem{acknowledgements}[theorem]{Acknowledgements}

\title{Level sets estimation and Vorob'ev expectation of random compact sets}

\author{Philippe Heinrich and Radu Stefan Stoica  and Viet Chi Tran\footnote{Laboratoire Paul Painlev\'{e} UMR CNRS 8524, Universit\'{e} Lille 1, 59 655 Villeneuve d'Ascq Cedex.}}

\begin{document}

\maketitle

\begin{abstract}
  The issue of a ``mean shape'' of a random set $X$ often arises, in
  particular in image analysis and pattern detection. There is no canonical definition but
  one possible approach is the so-called Vorob'ev expectation $\E_V(X)$, which is closely linked to quantile sets. In this paper,
  we propose a consistent and ready to use estimator of $\E_V(X)$ built from
  independent copies of $X$ with spatial discretization.  The control
  of discretization errors is handled with a mild regularity
  assumption on the boundary of $X$: a not too large `box counting'
  dimension. Some examples are developed and an application to cosmological data is presented.  
\end{abstract}

\noindent \textbf{keywords:} Stochastic geometry ; Random closed sets ; Level sets ; Vorob'ev
    expectation\\
\textbf{AMS codes:} Primary 60D05 ; Secondary 60F15 ; 28A80

\section{Introduction and background}
\label{intro}

The present paper proposes a ready to use estimator based on level sets, and that allows to tackle the question of approximating
the mean shape or expectation of a random set. The practical motivation of this work is given by pattern
recognition applications coming from domains like astronomy,
epidemiology and image
processing~\cite{StoiMartSaar07,StoiMartSaar10,StoiGayKret07,StoiDescZeru04}.
Several ways to define the expectation of a random set have been
developed in the literature (see \eg \cite{Molc99}). Our investigation
leads us to the Vorob'ev expectation, closely related to quantiles and
level-sets.\\
Let $(\Omega,\mathcal{A},\P)$ be a probability
space. Consider a random compact set $X$ in $[0,1]^d$ as a map from
$\Om$ to the class $\mathcal{C}$ of compact sets in $[0,1]^d$
measurable in the following sense (see \eg \cite{Molc99}):
$$\forall C\in\mathcal{C},\quad\{\om:X(\om)\cap C\ne \emptyset\}\in \mathcal{A}.$$

Robbins' formula, a straightforward consequence of Fubini's theorem, states that
$$\vm{X}=\int_\cube \P(x\in X)\l(dx),$$
where $\l$ is the Lebesgue measure.  To define a ``mean shape'', it is
thus natural to consider the coverage function $p(x)=\P(x\in X)$. For
$\alpha\in [0,1]$, the (deterministic) $\alpha$-level set of $p(x)$
is:
$$Q_\alpha=\{x\in [0,1]^d : p(x)>\alpha\}$$
or $\{p>\alpha\}$ for short. Choosing $\alpha$ such
that the volume of $Q_\alpha$ matches the mean volume of $X$ provides
the Vorob'ev expectation $\E_V(X)$ (see \textit{e.g.}
\cite{kovyazin,Molc99,molchanov,stoyanstoyan}).
\begin{definition}
\label{Fa*}
 Set $F(\alpha)=\lambda(Q_\alpha)=\lambda\{p>\a\}$ and
$$ \a^*=\inf\{\a\in [0,1]:F(\a)\le \vm{X}\}.$$
The Vorob'ev expectation  of $X$ is a Borel set $\E_V(X)$ such that
$\lambda(\E_V(X))=\E\lambda(X)$ and
\begin{equation} \{p>\a^*\} \subset \E_V(X)\subset
\{p\ge\a^*\}.\label{def:Vorob'ev}\end{equation}
\end{definition}
Uniqueness of $\E_V(X)$ is ensured if $F$ is continuous at $\a^*$, and
in this case one can choose $\E_V(X)=\{p\ge\a^*\}$ which is compact
since $p(x)$ is upper semi-continuous.  It can be shown (\cite[Theorem
2.3. p.177]{molchanov}) that $\E_V(X)$
minimizes $B\mapsto \E\lambda(B\bigtriangleup X)$ under the constraint
$\lambda(B)=\E\lambda(X)$, where $B\bigtriangleup X$ is the symmetric difference of $B$ and $X$.\\

Despite their very natural definitions, neither the $Q_\alpha$'s nor
$\E_V(X)$ are tractable for applications. First, the coverage
probability $p(x)$ is not always available in an analytical closed
form. Second, level sets can not be computed for all the points
$x\in [0,1]^d$ and discretization should be considered. \\
The first point has been tackled by a large literature: plug-in estimators obtained by replacing $p(x)$ with an empirical counterpart have been considered by Molchanov \cite{Molc87,Molc90,Molc98}, Cuevas \etal \cite{CuevFrai97,CuevGonzRodr06}. In particular, the latter establishes consistency of plug-in estimators for the $L^1$-norm under weak assumptions. In particular, the function $p(x)$ needs not being continuous. \\

The aim of this paper is to derive a consistent and implementable
estimator $K_{n,r}$, based on $n$ i.i.d. copies $X_1,\dots,X_n$ of $X$
and spatial discretizations with a grid of mesh $r$. We extend proofs
of consistency of plug-in estimators to estimators including a grid
discretization, and generalize existing works on level-sets for the
Vorob'ev expectation, when the level $\alpha$ is replaced by a level
$\alpha^*$ that depends on $X$. We prove strong consistency for the
symmetric difference and provide also convergence rates.

\section{Estimation of level sets}

\subsection{Plug-in estimation of level sets}

To define the plug-in estimators for sets $Q_\alpha$, let us consider the empirical
counterparts of $p(x)$:
$$    p_n(x)=\frac{1}{n}\sum_{i=1}^n\mathbf{1}_{\{x\in X_i\}}.$$
Then, the plug-in estimator is
$$Q_{n,\alpha}=\{p_n>\alpha\}.$$
In the literature, several distances can be used for closed
sets. Here, we are interested in the pseudo-distance defined for any
$A, B \in\mathcal{B}(\cube)$ by $\d(A,B)=\l(A\triangle B)$. Other
possibilities are the Hausdorff distance, for which similar
consistency results can be obtained, under additional assumptions on
the regularity of $p(x)$ (see \eg \cite{Molc90} and \cite[Th. 1 and
2]{CuevGonzRodr06}).  In \cite{Molc87,Molc90} for instance, $p(x)$ is
assumed continuous. The case where $p(x)$ is a density function has
also been much investigated. Kernel estimators may be used in the
plug-in estimation. Rates of convergence
\cite{Molc98,RigoVert09,Tsyb97} or asymptotic normality
\cite{masonpolonik} are considered with regularity assumptions on
$p(x)$. Adaptive estimation can also be performed
\cite{Singhscottnowak}.
\\

$L^1$-consistency is of practical interest and we follow the approach
developed for instance in \cite{CuevGonzRodr06}, and generalize their
results with grid approximations in view of numerical
implementation. For applications, we also aim at weak assumptions, in
particular, the continuity of the coverage function $p(x)$ should not
be required. As a result, the function $F$ is c\`{a}dl\`{a}g with possible
plateaus (constant regions). Note that $F$ is the survival function of the random
variable $p(U)$ with $U$ uniformly distributed on $\cube$.

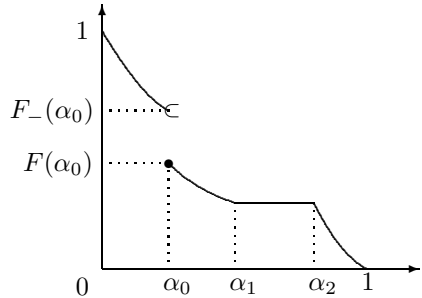
\begin{figure}[htpb]
\hspace{0.5cm}    \begin{picture}(150,110)(-20,-10)
      \put(0,0){\vector(1,0){120}} \put(0,0){\vector(0,1){100}}
      \put(-10,-10){0}
      \put(-10,87){1}
      \qbezier(0,90)(15,65)(25,60)
      \put(23,58){\footnotesize{$\subset$}}
      \dottedline{3}(0,60)(25,60)
      \put(-35,57){$F_-(\alpha_0)$}
      \dottedline{3}(0,40)(25,40)
      \put(-30,37){$F(\alpha_0)$}
      \put(23,38){\footnotesize{$\bullet$}}
      \qbezier(25,40)(35,30)(50,25)
      \dottedline{3}(25,0)(25,40)
      \put(23,-8){$\alpha_0$}
      \dottedline{3}(50,0)(50,25)
      \put(48,-8){$\alpha_1$}
      \put(50,25){\line(1,0){30}}
      \dottedline{3}(80,0)(80,25)
      \put(78,-8){$\alpha_2$}
      \qbezier(80,25)(90,05)(100,0)
      \put(98,-8){$1$}
    \end{picture}
\caption{Behavior of function $F(\alpha)=\lambda(Q_\alpha)$}
\label{figF}
\end{figure}
As an example, Figure~\ref{figF} shows how the function $F$ may
look like. It can be observed that the plateaus of $p(x)$ make the
discontinuities of $F$ while the discontinuities of $p(x)$ provide the
plateaus of $F$. The jump of $F$ at $\alpha_0$ corresponds to
the Lebesgue measure of the points with coverage probability equal to $\alpha_0$
$$
F_-(\alpha_0)-F(\alpha_0)=\lambda\{p=\alpha_0\}.
$$
The plateau $[\alpha_1,\alpha_2]$ of $F$ means that the set of points with
coverage probability between $\alpha_1$ and $\alpha_2$ is $\l$-negligible and
$\lambda\{p>\alpha_1\}=\lambda\{p>\alpha\}$ for all $\alpha\in
[\alpha_1,\alpha_2]$.

\subsection{Grid approximation of a level set}

The plug-in estimators $Q_{n,\alpha}$ are not satisfactory in view of applications, because they require
computation of $p_n(x)$ for all $x$. In all what follows, $r$ is a mesh introduced to work in $\cube\cap
r\Z^d$ rather than in $\cube$. We seek for a candidate $Q^r_{n,\alpha}$ close to $Q_{n,\alpha}$
for small $r$, but less greedy in computations.

\begin{definition}
\label{br}
For any Borel set $B$ in $\cube$ and $r\in 2^{-\N}$,  we call ``grid
approximation of $B$'' the set
$$  B^r=\bigsqcup_{x\in B\cap r\Z^d}[x,x+r)^d.$$
\end{definition}
As expected, some regularity of the border $\partial B$ has to be introduced (see \cite[p.38-39]{falconer}):
\begin{definition}
\label{dimb}
 Let
$$N_r(\partial B)=\mathrm{Card}\{x\in r\Z^d:[x,x+r)^d\cap \partial B\ne \emptyset \}.$$  The ``upper box counting dimension'' of $\partial B$ is
$$\dimb(\partial B)=\limsup_{r\to 0}\frac{\log N_r(\partial B)}{-\log r}.$$
\end{definition}
A straightforward consequence of these definitions is:
\begin{proposition}
\label{p1}
Assume that $\dimb(\partial B)<d$. For all $\eps >0$, there exists
  $r_\eps$ such that
$$0<r<r_\eps\Longrightarrow \d(B^r,B)\le r^{d-\dimb(\partial B)-\eps}.$$
\end{proposition}

The following result extends the result of Cuevas \etal \cite[Th. 3]{CuevGonzRodr06} by adding a discretization scheme:
\begin{proposition}\label{prop_approx_quantiles}
  Assume that $\dimb(\partial X)\leq d-\kappa$ with probability one for some $\kappa>0$. For all $\alpha$ such that $\lambda\{p=\alpha\}=0$,
  \begin{enumerate}
  \item[(i)] with probability $1$,
$$\lim_{\substack{r\to 0\\n\to \infty}}\d(Q_{n,\a}^r, Q_\alpha )=0.$$
\item[(ii)] for all $\eps >0$,
$$\E\,\d\big(Q_{n,\a}^r, Q_\alpha \big)\leq r^\kappa+
2\e^{-2n\eps^2}+F(\alpha-\eps)- F(\alpha+\eps).$$
  \end{enumerate}
\end{proposition}

\section{Estimation of $\E_V(X)$: Kovyazin's mean}
Let us consider the volume $F_n(\alpha) = \lambda\{p_n>\alpha\}$. Following the approach of \cite{kovyazin}, we
introduce first the empirical volume
$$\Lambda_n  = \frac{1}{n}\sum_{i=1}^n\lambda(X_i),$$
and next
$$\alpha^*_n  =  \inf\{\alpha\in [0,1]:F_n(\alpha)\le \Lambda_n \}.$$
\begin{definition}\label{ve}
The ``Kovyazin's mean'' is a Borel set $K_n$ defined by an empirical version of \eqref{def:Vorob'ev}:
$$\l(K_n)=\L_n \et \{p_n>\a_n^*\} \subset K_n\subset\{p_n\ge\a_n^*\}.$$
\end{definition}
The following revisits the approach of Kovyazin
\cite{kovyazin}:
\begin{theorem}
\label{th1}
Assume that $\l\{p=\a^*\}=0$. Then, with probability one,
$$\d(K_n,\E_v(X))\xrightarrow[n\to\infty]{}0. $$
\end{theorem}

\subsection{A grid approximation of $K_n$}
Let us consider the grid $\cube\cap
r\Z^d$. A first idea could be to take the grid
approximation of $B=X$ or $B=K_n$, but this would provide a set $B^r$ with mean volume
not necessarily equal to $\L_n$. Instead, we set
\begin{equation}
  \label{anr}
  \a_{n,r}^*=\inf\{\a\in [0,1]:\l\left(\{p_n>\a\}^r\right)\le \L_n\}
\end{equation}
and consider  a Borel set $K_{n,r}$  of volume $\L_n$ such that
\begin{equation}
\label{knr}
  \{p_n>\a_{n,r}^*\}^r\subset K_{n,r}\subset\{p_n\ge\a_{n,r}^*\}^r.
\end{equation}
The set $K_{n,r}$ is the implementable estimator that we propose. Its
definition amounts in the following procedure. First, we approximate
the expected volume of the random set $X$ by $\Lambda_n$, which gives
the number of cells to select. The latter are chosen according to
their estimated coverage, which is given by $p_n(x)$ for the cell
$[x,x+r)^d$ with $x\in r\Z \cap [0,1]^d$. The consistency and
convergence rate are provided by the following results.

\subsection{Consistency of $ K_{n,r}$}
Recall the definition of $\a^*$ in Definition~\ref{Fa*} and set
\begin{equation}
  \label{b*}
  \b^*=\sup\{\a\in [0,1]:F(\a)\ge \vm{X}\}.
\end{equation}
Our main theorem states that:
\begin{theorem}
\label{th2}
 Assume that
 \begin{enumerate}
 \item[(i)] $\P(\dimb(\partial X)\le d-\kappa)=1$ for some $\kappa >0$.
 \item[(ii)] The Lebesgue measures of $\{p=\a^*\}$ and $\{p=\b^*\}$ are zero.
 \end{enumerate}
 With probability one,
$$\d(K_{n,r},\E_V(X))\xrightarrow[\substack{n\to \infty\\ r\to 0}]{}0.$$
\end{theorem}
For the proof, we write that
$$\d(K_{n,r},\E_V(X))\le  \d(K_{n,r},K_n)+\d(K_n,\E_V(X))$$
and use Theorem~\ref{th1} and the two following lemmas to conclude.

\begin{lemma}
\label{lm1}
Assume that $\P(\dimb(\partial X)\le d-\kappa)=1$ for some $\kappa >0$. Then, with probability one
$$\a^*\le\liminf_{\substack{n\to \infty\\ r\to 0}} \a_{n,r}^* \le
  \limsup_{\substack{n\to \infty\\ r\to 0}} \a_{n,r}^*\le \b^*,$$
and as a sort of particular case, we also have
$$\a^*\le\liminf_{n\to \infty} \a_n^* \le
  \limsup_{n\to \infty} \a_{n}^*\le \b^*.$$
\end{lemma}
\begin{lemma}
\label{lm2}
  Assume $\P(\dimb(\partial X)\le d-\kappa)=1$ for some $\kappa >0$. With
  probability one,
$$\limsup_{\substack{n\to \infty\\ r\to 0}}\d(K_{n,r},K_n)\le
2\Big[\lim_{\substack{\a\to \a^*\\ \a<\a^*}}F(\a)-F(\b^*)\Big].$$
\end{lemma}

\section{Some examples}\label{section4}

\subsection{The region covered by a boolean model}
As an example of multi-dimensional random set, let us consider the
Boolean model~\cite{lantuejoul,Molc97,StoyKendMeck95}, constructed as follows. First, take a Poisson point process
$\Pi_\mu$ in $\R^{d}$ ; the parameter $\mu$ is a locally finite measure on
$\R^d$ called intensity measure. Next consider a sequence of i.i.d compact sets $(\Xi_x)_{x\in \R^d}$, which is independent of $\Pi_\mu$. Finally, replace each point $x$ of
$\Pi_\mu$ by the shifted corresponding set $x+\Xi_x$. The
resulting union set
$$\Xi= \bigcup_{x \in \Pi_\mu}(x + \Xi_x)$$
is the Boolean model. The points $x$ are called germs and the random set $\Xi_0$ is the
`typical' grain of the model.  The Boolean model is also called the Poisson germ-grain
model. To avoid trivial models where $\Xi=\R^d$ almost surely, the $d$-th
moment of the radius of the circumscribed circle of $\Xi_0$ must be
finite.

The distribution of the Boolean model is uniquely determined by its
capacity functional (see~\cite{Molc97,StoyKendMeck95})
\begin{equation}
T_{\Xi}(C)= \P(\Xi \cap C  \neq \emptyset) = 1 - \exp\left[-\E\mu(\breve{\Xi}_0\oplus C)\right],
\label{capacityFunctional}
\end{equation}
where $C$ ranges over compact sets in $\R^d$, $\breve{\Xi}_0$ is the
symetric of typical grain with respect to the origin and $\oplus$ the
Minkowski addition.

Assume that the Boolean model is observed through the window $\cube$ so that we focus on the random compact set
$$X=\Xi\cap\cube,$$
that is the covered region in the cube by $\Xi$.
That $X$ is indeed compact requires some (known) work. Briefly, if $\Xi$ is viewed  as a random counting measure
on the class $\mathcal{C}$ of compact sets in $\R^d$, its intensity measure
will be locally finite if and only if
$\E\mu(\breve{\Xi}_0\oplus C)<\infty $ for every compact $C$. In these conditions, only finitely
many sets  $x+\Xi_x$ will hit the cube $\cube$ almost surely (see
\cite[sections 3 and 4]{SchWeil}).
\subsubsection{Stationary case}
\label{exemple_homogeneous_case}
This type of Boolean model refers to a translation-invariant intensity
measure $\mu$ on $\R^d$ so that $\mu(dx) = m \l(dx)$ for some $m
>0$. In the following, the considered grains $\Xi_x$ are balls
of random radius $R$ such that $\E(R^d) < \infty$. The coverage
probability $p(x)$ of $X$ is obtained from~\eqref{capacityFunctional} by
taking $C=\{x\}$ with $x\in \cube$:
$$ p(x) = T_{\Xi}(\{x\})=1 - \exp\left[ - m \l(B(0,1))\E( R^d)\right].
$$
Set for short $c_{m,d}=1 - \exp\left[ - m \l(B(0,1))\E( R^d)\right]$. Since $p(x)$ is the constant $c_{m,d}$, it's easy to conclude that
$$
 F(\alpha)=\ind_{[0,c_{m,d})}(\alpha)\quad\text{and}\quad\a^*=c_{m,d}.
$$
Thus, $Q_{\a^*}=\emptyset$ while $Q_{\a^*-\varepsilon}=[0,1]^d$ for any $\varepsilon>0$. In this
case, $\E_v(X)$ is not unique and any measurable set with volume
$\E\lambda(X)$ is a possible Vorob'ev expectation. Also, the
assumption $\lambda\{p=\a^*\}=0$ does not hold here and neither
does Theorem~\ref{th1}.

\subsubsection{Non-stationary case}\label{section:nonstationary}
This type of Boolean model may be obtained for locally finite intensity
measures $\mu(dx) = m(x)\lambda(dx)$ with a non-constant positive function $m(x)$. We assume that the radius probability distribution $\P\circ R^{-1}$
admits a continous density $g$ with respect to the Lebesgue measure such that
\begin{equation}
\forall x\in\cube,\quad \int_\cube g(|x-y|)\frac{x-y}{|x-y|}\mu(dy) \neq 0.
\label{eqH4}
\end{equation}
The coverage probability of $X$ requires the computation of the function
\begin{eqnarray*}
\phi(x)  &=& \E\mu\left(\breve{\Xi}_0\oplus \{x\}\right)\\
 &=& \E\mu\left(x+B(0,R)\right)\\
& =& \int_{[0,1]^d}\P(R>|x-y|)
\mu(dy).
\end{eqnarray*}
For $\alpha\in (0,1)$, the level set $Q_\alpha$ consists of points $x$ satisfying:
\begin{align*}
p(x)>\alpha \quad  \Leftrightarrow \quad & \exp\big(-\phi(x)\big)<1-\alpha\\
\quad  \Leftrightarrow  \quad & \phi(x)>\ln\Big(\frac{1}{1-\alpha}\Big).
\end{align*}
Computing the gradient $\nabla_x \P(R>|x-y|)$ for $x\ne y$ gives the integrand of~\eqref{eqH4}. Since $g$ is continuous, $\phi$ is of
class $\Co^1$ with the non-vanishing derivative \eqref{eqH4}. Moreover, the boundary $\partial Q_\alpha$ coincide with
the set $\{p=\alpha\}$ and is obtained as the solution of
$$\phi(x)=\ln\Big(\frac{1}{1-\alpha}\Big).$$
Using the implicit function theorem, the boundary of $Q_\alpha$ is a
$\Co^1$-manifold of dimension $d-1$: it can be locally parameterized
by $\Co^1$-functions of $d-1$ variables. Assumption $\lambda\{p=\alpha\}=0$
 is satisfied. In this case,
Proposition~\ref{prop_approx_quantiles} holds for any $\alpha\in
(0,1)$ and Theorems~\ref{th1} and \ref{th2} apply, which yields consistency and a rate of convergence for the estimator $K_{n,r}$.

\subsubsection{Boolean model with atoms}

In this section, we provide an example where the coverage function exhibits discontinuities. Consider the model $\Xi$ of Section \ref{section:nonstationary} with in addition an independent random vector $(U_1,\ldots, U_N)$ of independent Bernoulli random variables. Fix points $y_1,\ldots, y_N$ in $\cube$ and $r_0>0$ and set
$$\Upsilon=\Xi \cup \bigcup_{i: U_i=1}B(y_i,r_0).$$
We have thus added a random union of deterministic balls to the model $\Xi$. In the case where some of the $\P(U_i=1)$ are equal to 1, this model is reminiscent of the conditional Boolean model (see Lantu\'{e}joul \cite{lantuejoul}).\\
The coverage probability of $X=\Upsilon\cap \cube$ is:
$$p(x)=1-e^{-\phi(x)}\prod_{i :  B(y_i,r_0)\ni x}\P(U_i=0).$$
 and $p(x)>\alpha$ is equivalent to
\begin{equation}
\phi(x)-\sum_{i :  B(y_i,r_0)\ni x}\ln\P(U_i=0)>\ln\Big(\frac{1}{1-\alpha}\Big).\label{etape}
\end{equation}
Set $\mathcal{D}=\cup_{i=1}^N\partial B(y_i,r_0)$. If we work on each connected component of $[0,1]^d\setminus\mathcal{D} $, the sum in the l.h.s. of
\eqref{etape} is constant and we are reduced to the case considered
in Section \ref{section:nonstationary}. Besides, $p(x)$ is discontinuous on the set $\mathcal{D}$. Thus $F(\alpha)$ exhibits plateaus. If $p(x)$ is continuous at $\alpha^*$, then the Vorob'ev expectation is
unique and consistency and rate of convergence can be obtained for
$K_{n,r}$. Else $p(x)$ is discontinuous at $\alpha^*$ and there is no
uniqueness of the Vorob'ev expectation.

\subsection{Example on real data}
\begin{figure}[!htbp]
\centering
\epsfxsize = 8cm \epsffile{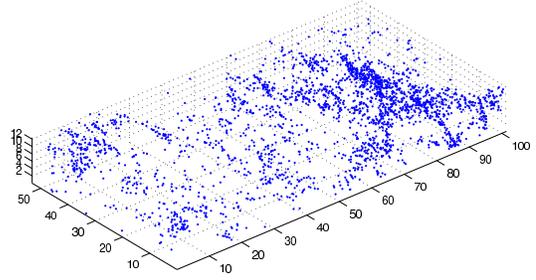}
\caption{NGP150 sample from the 2dFGRS survey.}
\label{galaxyCatalogue}
\end{figure}

The filamentary network is one of the most important features in
studying galaxy distribution~\cite{MartSaar02}. The point field in
Figure~\ref{galaxyCatalogue} is given by the galaxy centers of the
NGP150 sample obtained from the 2 degree Field Galaxy Redshift Survey
(2dFGRS)~\cite{CollEtAl01}. It can be easily noticed that the galaxy
positions are not spread uniformly in the observation window and the
filamentary structure is observable simply by eye investigation.

The automatic delineation of these structure is an important challenge
for the cosmological community. The authors in~\cite{StoiMartSaar07,StoiMartSaar10} models these filamentary structures by a cloud $X$ of interacting cylinders $C(x,\eta)$ which are encoded by their centers $x\in K$ compact subset of $\R^3$ and their directions $\eta\in \mathbb{S}^3_+$ the unit sphere with positive $z$-ordinate. To obtain such an object, these authors use a point process $\Gamma$ in $\R^3\times \mathbb{S}^3_+$ which has density
\begin{equation}
h(\gamma|\theta) = \frac{\exp[-U_d(\gamma|\theta)-U_i(\gamma |\theta)]}{\alpha(\theta)},
\label{bisousModel}
\end{equation}with respect to the Poisson point process $\Pi_\mu$ on $\R^3\times \mathbb{S}^3_+$ with the Lebesgue measure as intensity $\mu$.
The energy function of the model is given by the sum of two terms that depend on parameters $\theta\in \Theta$ a compact subset of $\R^m$.
The term $U_{d}(\gamma|\theta)$ depends on the observed galaxy fields and is called the data energy. It expresses that the point process giving the centers $x$ is spatially inhomogeneous: more points will be drawn in region where observations are dense. The term $U_{i}(\gamma|\theta)$ is called the
interaction energy and it is related to the connection and the
alignement of the cylinders forming the pattern.\\

The random set that is considered is hence:
$$X=X(\Gamma)= \Big(\bigcup_{(x,\eta)\in \Gamma}C(x,\eta)\Big)\cap K.$$
%In the sequel, we denote by $p(X|\theta)$ the distribution of $X$ given the parameters $\theta$, that is the image measure of $\Gamma$ through the application $\Gamma\mapsto X$.
For full complete details concerning the model and the simulation
dynamics, the reader should refer to
in~\cite{StoiMartSaar07,StoiMartSaar10}. To fit our previous developped
framework, the compact $K$ should be chosen as (or imbedded in)
$[0,1]^d$ by using homothetical scaling of the data. But this is completely inessential.

The model~\eqref{bisousModel} behaves as follows. If only the
interaction energy is used, the model simulates a connected network
which is independent of the data. If only the data term is used the
model locates the filamentary regions but the cylinders are clustering
in a disordered way and they neither connect nor form a network. The
data energy acts as a non-stationary Poisson marked point process.
The interaction energy acts as distance based pair-wise interaction
marked point process and it plays the role of a regularization term or
a prior for the filamentary network.\\

Heuristically, the idea in \cite{StoiMartSaar07,StoiMartSaar10} is to estimate the parameter $\theta$, which is unknown, and simulate realizations of $X$ from which information on the distribution of $X$ can be obtained. For instance, if $X_1,\dots,X_n$ are realizations of $X$, their Vorob'ev mean can help us identify the patterns of the filamentary network by providing a mean shape appearing under the distribution of $X$. Tunning a model as \eqref{bisousModel} is not
always a very simple task. Stoica et al. use a Bayesian frameworks and introduce a prior $\varrho(\theta)$ for the model parameters. Using a Metropolis-Hastings algorithm, throwing the burn-off period away and choosing simulations that are sufficiently distant, they generate a non-correlated sequence $(\theta_i,X_i)_{i\in \{1,\dots n\}}$ that should be identically distributed if the algorithm converged correctly. The Metropolis-Hastings algorithm uses simulated annealing and depends on a temperature parameter $T$, and several strategies can be applied.\\
Once this is done, notice that the coverage function writes as:
\begin{equation}
p(x)  =  \E\left(\int_{\Theta} \1_{\{x \in X(\Pi_\mu) \}} h(\Pi_\mu|\theta)\varrho(\theta)d\theta\right).
\label{coveragePattern}
\end{equation}
The Lebesgue measure of the level sets induced
by~\eqref{coveragePattern} is given by the function
\begin{eqnarray*}
F(\alpha) %& = & \lambda\{Q_\alpha\}=\lambda\{x \in K : p(x) > \alpha\} \nonumber \\
& = & \int_{K}\1_{\{p(x) > \alpha\}}\lambda(dx).
\label{functionPattern}
\end{eqnarray*}
Since, $p(x)$ is not available in closed analytical form it is not a
trivial task to check the continuity of $F(\alpha)$ to check if we are in the assumptions of Theorem \ref{th2}.
%Therefore, two
%simulation tests were performed for the data set in
%Figure~\ref{galaxyCatalogue} in order to verify the availability of
%the continuity hypothesis of $F(\alpha)$. The model and the dynamics
%parameters were set as in~\cite{StoiMartSaar07}.

\subsubsection{Fixed temperature}

If the temperature remains constant during the Metropolis-Hastings algorithm, then the algorithm is expected to produce a sequence of non-correlated random variables $(\theta_i,\Gamma_i)$ with density $h(\gamma|\theta)\varrho(\theta)$. The sequence $(X_i)_{i\in \{1,\dots,n\}}$ is then considered to be identically distributed.\\

\begin{figure}[!htbp]
\centering
\epsfxsize = 8cm \epsffile{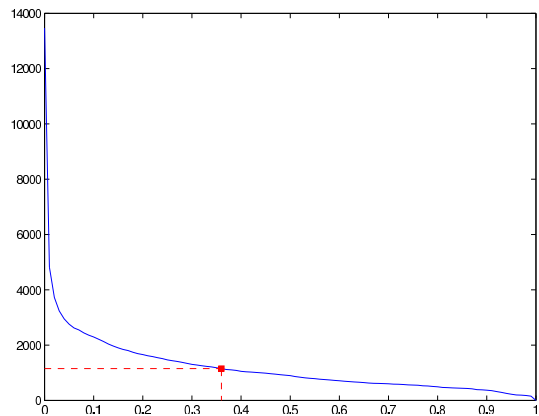}
\caption{Monte-Carlo computation of the function $F(\alpha)$ for the
  galaxy catalogue NGP150 using a simulation dynamics at fixed
  temperature.}
\label{functionTFixed}
\end{figure}

The first empirical test consists in running the simulation
dynamics at fixed temperature $T=1$, and we simulate a sequence of $n=1000$ realizations of $X$. The Monte Carlo counter-part of $F(\alpha)$ is
shown in Figure~\ref{functionTFixed}. No observable discontinuity can
be detected by simple visual inspection. Hence, using the same samples
we obtain the empirical values $\Lambda_n = 1158$ and
$\alpha_n^{*} = 0.36$.

\begin{figure}[!htbp]
\centering
\epsfxsize = 8cm \epsffile{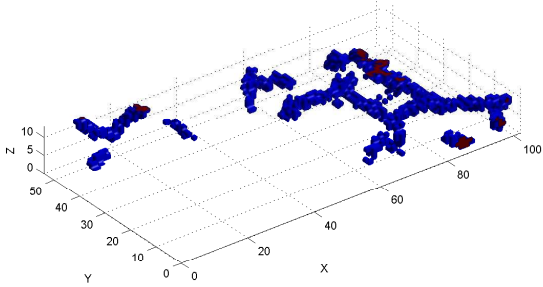}
\caption{Estimator of the Vorob'ev expectation for the
  galaxy catalogue NGP150 using a simulation dynamics at fixed
  temperature.}
\label{vorobevTFixed}
\end{figure}

 The estimator of the Vorob'ev expectation is shown in
 Figure~\ref{vorobevTFixed}. Even if running the simulation dynamics
 at fixed temperature does not provide an estimator of the
 filamentary pattern, the estimator of the Vorob'ev expectation can be
 used to assess the presence of a filamentary pattern in the data under
 the assumption that the model parameters are correct. The empirical
 computation of the sufficient statistics of the model can be also
 used in order to validate the presence of the filamentary pattern and
 also to give a more detailed morphological
 description~\cite{StoiGayKret07,StoiMartSaar07,StoiMartSaar10}.

\subsubsection{Decreasing temperature}

We now start with a low temperature $T_0$ and choose the cooling schedule along the Metropolis-Hastings algorithm as
\begin{equation}
T_n = \frac{T_0}{\log n + 1}.
\label{coolingSchedule}
\end{equation}
Under these circumstances, the sequence $(\theta_i,\Gamma_i)_{i\in \{1,\dots,n\}}$ can be considered as non-correlated random variables that are uniformly distributed on the set  maximizing the joint probability density or minimizing
the total energy of the system (see \cite{StoiGregMate05})
\begin{eqnarray}
(\widehat{\theta},\widehat{\gamma}) & = & \arg\max_{\Omega\times \Theta}
  h(\gamma|\theta)\varrho(\theta) \nonumber \\
& = &\arg\min_{\Omega \times \Theta}\left\{
  \frac{U_{\dd}(\gamma|\theta)+U_{i}(\gamma|\theta)}{\alpha(\theta)} - \log
  \varrho(\theta)\right\},
\label{estimatorPattern}
\end{eqnarray}
where $\Omega$ denotes the set of simple point measures on $\R^3\times \mathbb{S}^3_+$. Clearly, the solution we obtain is not unique. \\

\begin{figure}[!htbp]
\centering
\epsfxsize = 8cm \epsffile{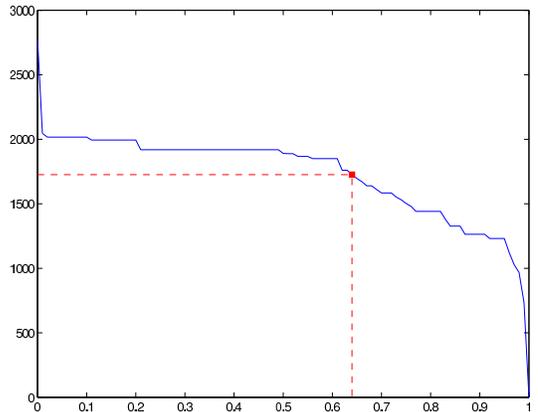}
\caption{Computation of the function $F(\alpha)$ for the
  galaxy catalogue NGP150 using a simulated annealing algorithm.}
\label{functionSAnnealing}
\end{figure}

The second test was to obtain $n=1000$ samples using the simulated
annealing algorithm with $T_0=1$.  The approximated $F(\alpha)$
function is shown in Figure~\ref{functionSAnnealing}. Comparing with
the previous case, the function $F(\alpha)$ looks less smooth. Two
other approximations are done. First, the simulated annealing
simulates a non-homogeneous Markov chain. Second, the cooling schedule
fixes the sequence of distributions converging towards the uniform
distribution on the configuration sub-space maximizing
$h(\gamma|\theta)\varrho(\theta)$. Again using the same samples we
obtain the empirical values $\Lambda_n = 1730.6$ and $\alpha_n^{*} =0.64$.

\begin{figure}[!htbp]
\centering
\epsfxsize = 8cm \epsffile{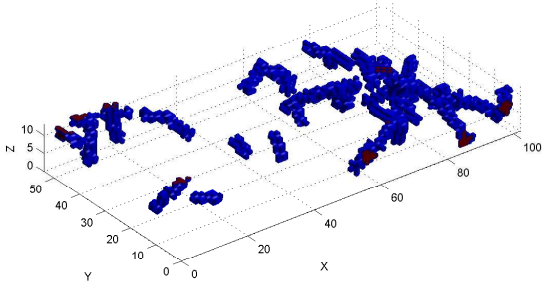}
\caption{Estimator of the Vorob'ev expectation for the
  galaxy catalogue NGP150 using a simulated annealing algorithm.}
\label{vorobevSAnnealing}
\end{figure}

The estimator of the Vorob'ev expectation computed using the simulated
annealing algorithm is shown in
Figure~\ref{vorobevSAnnealing}. Comparing with the first experience,
the filamentary pattern has a bigger volume. As a consequence, this method is more useful for detecting average patterns
of the filaments, while one would prefer the method with fixed temperature for estimation.

\bigskip The presented applications lead towards new questions
concerning the level sets estimators. The first question consists in
determining a class of models that allow computation of level sets and
Vorob'ev expectation in the setting of this paper. A very intuitive
answer leads to the following recommendations of using continuous
intensity functions and interaction potentials. Clearly, this may be a
quite restrictive condition. Hence, the second question would be
whether continuous priors on the model parameter guarantee an
appropriate behavior of the $F(\alpha)$ function. And finally a third
question is whether there exists an optimal cooling schedule
guaranteeing the working hypotheses for the construction of level sets
based estimators. From our experiments, it seems for instance that
constant temperature is good for estimating parameters while simulated
annealing is better to detect patterns.

\appendix
\section{Appendix: Proofs}
Several times, the Strong Law of Large numbers in the separable Banach
space $L^1(\cube)$ ($L^1$ SLLN) will be used (see \eg
\cite[chap. 7]{LedTal}).
\subsection*{Proof of Proposition~\ref{p1}}
It is not difficult to see that $B\bigtriangleup B^r$ is included in the union of
cells $[x,x+r)^d$ ($x\in r\Z^d\cap\cube$) that meet $\partial B$. As a result,
we have
$$\l(B\bigtriangleup B^r)\le N_r(\partial B)r^d.$$
Besides, given $\eps>0$, we get from definition \eqref{dimb} that for $r$
small enough,
$$ N_r(\partial B)\le r^{-\dimb(\partial B)-\eps}.$$
These two inequalities provide the result.

\subsection*{Proof of Lemma~\ref{lm1}}
Consider $\a<\a^* $ such that $F$ is continuous at $\a$. We show that
$$\P\Big( \a\le \liminf_{\substack{n\to \infty\\ r\to 0}}
\a_{n,r}^*\Big )=1.$$
Otherwise, on a set of positive probability, we have
$\a> \liminf_{\substack{n\to \infty\\ r\to 0}}\a_{n,r}^*$ and there
exists a sequence $(n_k,r_k)_{k\ge 1}$ with $n_k\to\infty$ and $r_k\to 0$ such
that $\a>\a_{n_k,r_k}^*$. Then, from  the definition \eqref{knr},
\begin{eqnarray*}
  \L_{n_k} &=&\l(K_{n_k,r_k})\\
 &\ge & \l\left(\{p_{n_k}>\a\}^{r_k}\right)\\
&\ge &
\l\left(\{p_{n_k}>\a\}\right)-\d\left(\{p_{n_k}>\a\}^{r_k},\{p_{n_k}>\a\}\right)\\
&\ge &
\l\left(\{p_{n_k}>\a\}\right)-\sup_n\d\left(\{p_n>\a\}^{r_k},\{p_n>\a\}\right).
\end{eqnarray*}
Note that $\partial\{p_n>\a\}\subset \bigcup_{i=1}^n\partial
X_i$ since $p_n$ is locally constant on the complementary of
$\bigcup_{i=1}^n\partial X_i$. Moreover, since the ``upper box-counting dimension'' $\dimb(\cdot)$ has
monotonic and stability properties (see \cite{falconer}), we get whatever $\a$, for the constant $\kappa$ introduced in the assumptions,
\begin{eqnarray*}
\dimb(\partial\{p_n>\a\}) &\le& \dimb\Big(\bigcup_{i=1}^n\partial
X_i\Big)\\
&=&\max_{1\le i\le n}\dimb\Big(\partial
X_i\Big)\le d-\kappa.
\end{eqnarray*}
It follows then from Proposition~\ref{p1} that for $k$ large enough,
\begin{equation}
  \label{rkappa}
  \sup_n\d\left(\{p_n>\a\}^{r_k},\{p_n>\a\}\right)\le r_k^{\kappa/2},
\end{equation}
and thus
$$\L_{n_k}\ge \l\{p_{n_k}>\a\}-r_k^{\kappa/2}.$$
Taking the limit as $k\to \infty$ and invoking the $L^1$-SLLN for
$(p_n)_{n\geq 1}$ , we obtain
\begin{equation}
\vm{X}\ge \l\{p>\a\}=F(\a).\label{but1}
\end{equation}
This contradicts the definition of $\a^*$. As a consequence, we have
$$\P\Big( \a^*\le \liminf_{\substack{n\to \infty\\ r\to 0}}
\a_{n,r}^*\Big )=1.$$
Very similarly, we can prove that
$$\P\Big(\limsup_{\substack{n\to \infty\\ r\to 0}}
\a_{n,r}^*\le \b^*\Big )=1.$$
The particular case without the mesh $r$ is in the same spirit: we can show that there exists a sequence $n_k\rightarrow\infty$ such that $\Lambda_{n_k}\geq \lambda(\{p_{n_k}>\alpha\})$ which leads to \eqref{but1}.

\subsection*{Proof of Lemma~\ref{lm2}}
From $\d(A,B)=\l(A)-\l(B)+2\l(B\setminus A)$, and since $K_{n,r}$ and $K_n$ have the same volume, we deduce
\begin{eqnarray*}
 \frac{1}{2}\d(K_{n,r},K_n) &  \le & \l\left(\{p_n\ge \a_{n,r}^*\}^r\setminus
    \{p_n> \a_{n}^*\}\right)\\
&\le & \l\left(\{p_n\ge \a_{n,r}^*\}^r\setminus
    \{p_n\ge \a_{n,r}^*\}\right)\\
& & \phantom{\{p_n\ge \a\}}+\l\left(\{p_n\geq \a_{n,r}^*\}\setminus
    \{p_n> \a_n^*\}\right)\\
&\le &\d\left(\{p_n\ge \a_{n,r}^*\}^r,\{p_n\ge
  \a_{n,r}^*\}\right)\\
& & \phantom{\{p_n\ge \a\}ccccccccc}+\l\{\a_{n,r}^*\le p_n\le \a_n^*\}.
\end{eqnarray*}
The first term in the majoration is with probability one less than
$r^{\kappa/2}$ for $r$ small enough, as in \eqref{rkappa}. For the second one,
let $\a <\a^*$  and $\b>\b^*$ ; by Lemma \ref{lm1}, there exists $n_{\a,\b}\ge 1$ and $r_{\a,\b} >0$ such that for all $n\ge n_{\a,\b}$ and $r\in (0,r_{\a,\b})$,
$$\l\{\a_{n,r}^*\le p_n\le \a_n^*\}\le \l\{\a< p_n\le\b\},$$
and again by the $L^1$-SLLN for $(p_n)_{n\geq 1}$,
$$\l\{\a < p_n\le\b\}\xrightarrow[n\to \infty]{}\l\{\a<p\le\b\}=F(\a)-F(\b),$$
provided $F$ is continuous at $\a$ and $\b$. We deduce that with probability one
$$\limsup_{\substack{n\to \infty\\ r\to 0}}\d(K_{n,r},K_n)\le 2\left[F(\a)-F(\b)\right],$$
and the proof is easily completed.

\subsection*{Proof of  Theorem~\ref{th1}}
Consider Borel sets $A,B$, $\a,\b\in[0,1]$ and coverage functions $p,q$ such
that
$$\{p>\a\}\subset A\subset \{p\ge \a\}\et \{q>\b\}\subset B\subset \{q\ge
\b\}.$$
From $\d(A,B)=\l(A)-\l(B)+2\l(B\setminus A)$, we deduce
\begin{multline*}
\d(A,B)\le |\l(A)-\l(B)|\\ +2\left[\l\{p\ge \a,q\le\b\}\wedge\l\{p\le
  \a,q\ge\b\}\right],
\end{multline*}
and for $A=K_n$ and $B=\E_V(X)$, we get
\begin{equation}
  \label{maj}
  \d(K_n,\E_V(X))\le |\L_n-\l(\E_V(X))|+2m_n
\end{equation}
with
$$m_n=\l\{p_n\ge \a_n^*,p\le\a^*\}\wedge\l\{p_n\le
  \a_n^*,p\ge\a^*\}.$$
The first term in \eqref{maj}  converges a.s. towards $0$ by the SLLN. The
second one requires to distinguish two cases:
\begin{description}
\item[$\a^*=\b^*$:] It's enough to prove that, with probability one,
  \begin{equation}
    \label{cvl1}
  p_n-\a_n^*\xrightarrow[n\to\infty]{L^1(\cube)} p-\a^*.
  \end{equation}
Indeed together with the assumption $\l\{p=\a^*\}=0$ this yields
$$\l\{p_n- \a_n^*\le 0, p\ge\a^*\}\xrightarrow[]{}\l\{p-\a^*\le 0, p\ge\a^*\}=0,$$
and thus $m_n$ tends to $0$ as $n\to\infty$.
But \eqref{cvl1} stems from Lemma~\ref{lm1} that gives
$\a_n^*\xrightarrow[]{}\a^*$  with probability one  and
$p_n\xrightarrow[]{L^1(\cube)} p$ also with probability one.
\item[$\a^*<\b^*$:] We distinguish three subcases. If $\a^*_n<\a^*$, then
  \begin{eqnarray*}
    m_n &\le & \l\{p_n\le\a_n^*,p\ge\a^*\}\\
&\le & \l\{p_n\le\a_n^*,p\ge\b^*\}+\l\{\a^*\le p<\b^*\}\\
& = & \l\{p_n\le\a_n^*,p\ge\b^*\}\\
& & \phantom{ccccccccccccc}+\l\{p=\a^*\}+(F(\a^*)-F(\b^*-)).
 \end{eqnarray*}
The two last terms  are zero, thus:
   $$ m_n \leq  \l\{|p_n-p|\geq \b^*-\a^*\}
   \leq \frac{\| p_n-p\|_{L^1}}{\b^*-\a^*},$$
by using the Markov
   inequality. If $\a^*_n>\b^*$, then starting from $m_n\leq
   \l\{p_n\ge \a^*_n, p\le \a^*\}$ we obtain the same
   upper-bound. Finally, when $\a^*\le\a^*_n\le\b^*$,
  \begin{eqnarray*}
(\b^*-\a^*)m_n &\leq & (\b^*-\a^*_n)\l\{p_n\le \a^*_n,p\ge \b^*\}\\
&  &\phantom{cccccc} +(\a^*_n-\a^*)\l\{p_n\ge \a^*_n,p\le \a^*\}\\
  &\leq & \|p_n-p\|_{L^1}.
  \end{eqnarray*}
In the three cases, $m_n$ is bounded by $\|p_n-p\|_{L^1}/(\b^*-\a^*)$ which converges to zero a.s.
\end{description}

\subsection*{Proof of Proposition~\ref{prop_approx_quantiles}}

By triangular inequality and Proposition~\ref{p1},
$$
\E\,\d\big(Q^r_{n,\alpha}, Q_\alpha \big)\leq
r^{\kappa}+\E\,\d\big(Q_{n,\alpha},Q_\alpha\big).
$$
It remains to control the last term $ \E\,\d\big(Q_{n,\alpha},Q_\alpha\big)$. We have
\begin{eqnarray*}
\d(Q_{n,\alpha},Q_\alpha) &= &\lambda\big\{p_n> \alpha,\,p\le
\alpha\big\}+ \lambda\big\{p_n\le \alpha,\,p> \alpha\big\}\nonumber\\
&\leq &\lambda\big\{x \in [0,1]^d : |p_n(x)-p(x)|\geq |p(x)-\alpha|\big\}.
\end{eqnarray*}
Taking expectation and using Fubini's theorem and Bernstein's
inequality, one gets for all $\eps>0$
\begin{eqnarray*}
  \E\,\d\big(Q_{n,\alpha},Q_\alpha\big)&\le &\int_{[0,1]^d}\P\{|p_n(x)-p(x)|\geq
 |p(x)-\alpha|\}\lambda(dx) \nonumber\\
&\le &\int_{\{|\alpha-p|\ge \eps\}} 2\e^{-2n|\alpha-p(x)|^2}\lambda(dx)\\
& &\phantom{cccccccccccccccccccc} +
\int_{\{|\alpha-p|< \eps\}}\lambda(dx)\nonumber\\
&\le & 2\e^{-2n\eps^2}+F(\alpha-\eps)-
F(\alpha+\eps).
\end{eqnarray*}

%\section{Section title}
%\label{sec:1}
%Text with citations \cite{RefB} and \cite{RefJ}.
%\subsection{Subsection title}
%\label{sec:2}
%as required. Don't forget to give each section
%and subsection a unique label (see Sect.~\ref{sec:1}).
%\paragraph{Paragraph headings} Use paragraph headings as needed.
%\begin{equation}
%\end{equation}

%% For one-column wide figures use
%\begin{figure}
%% Use the relevant command to insert your figure file.
%% For example, with the graphicx package use
%  \includegraphics{example.eps}
%\end{figure}
%%
%% For two-column wide figures use
%\begin{figure*}
%% Use the relevant command to insert your figure file.
%% For example, with the graphicx package use
%  \includegraphics[width=0.75\textwidth]{example.eps}
%\end{figure*}
%\begin{table}
%\begin{tabular}{lll}
%\hline\noalign{\smallskip}
%\end{tabular}
%\end{table}

\begin{acknowledgements}
The authors are very grateful to V. Martinez and E. Saar for providing the
astronomical data and to Ilya Molchanov and to an anonymous
referee for useful comments.The authors also thank the "G\'eom\'etrie stochastique" workgroup of the University Lille 1 for interesting discussions.
\end{acknowledgements}

{\footnotesize
\providecommand{\noopsort}[1]{}

}

\end{document}